\documentstyle{amsppt}
\magnification=1200
\hoffset=-0.5pc
\vsize=57.2truepc
\hsize=38truepc
\nologo
\spaceskip=.5em plus.25em minus.20em
\define\Gammat{\tau}
\define\Sigm{\Cal S}
\define\fra{\frak}
\define\Ri{R}

\define\Bobb{\Bbb}
\define\barakont{1}
\define\batviltw{2}
\define\batvilfo{3}
\define\batavilk{4}
\define\berikash{5}
\define\berikafi{6}
\define\berikafo{7}
\define\bogomone{8}
\define\bogomtwo{9}
\define\cartanon{10}
\define\cartantw{11}
\define\cartanei{12}
\define\chenfou{13}
\define\eilmactw{14}
\define\gerstenh{15}
\define\gerschtw{16}
\define\gugenhtw{17}
\define\gugenhth{18}
\define\gugenlam{19}
\define\gulasta{20}
\define\gulstatw{21}
\define\gugenmun{22}
\define\hainone{23}
\define\haintond{24}
\define\halpstas{25}
\define\halpsttw{26}
\define\hirschon{27}
\define\perturba{28}
\define\cohomolo{29}
\define\modpcoho{30}
\define\intecoho{31}
\define\bv{32}
\define\banach{33}
\define\huebkade{34}
\define\huebstas{35}
\define\kadeifou{36}
\define\kosmathr{37}
\define\kuranone{38}
\define\kurantwo{39}
\define\maninfiv{40}
\define\moorefiv{41}
\define\munkholm{42}
\define\quilltwo{43}
\define\sanebfou{44}
\define\sanebfiv{45}
\define\schlstas{46}
\define\stashnin{47}
\define\tianone{48}
\define\todortwo{49}
\define\wittetwe{50}

\topmatter
\title
Berikashvili's functor $D$
and the deformation equation
\endtitle
\dedicatory
Dedicated to Nodar Berikashvili
on the occasion of his seventieth birthday
\enddedicatory
\author
Johannes Huebschmann
\endauthor
\abstract
{Berikashvili's
functor $\Cal D$ 
defined in terms of twisting cochains 
is related to deformation theory,
gauge theory, Chen's formal power series connections,
and the master equation in physics.
The idea is advertised that
some unification
and understanding of the links between
these topics
is provided by the notion of twisting cochain
and the idea of classifying twisting cochains.}
\endabstract
\subjclass{
13D10
14B12
14J32
16S80  
17A36
17B55 
17B56 
17B65
17B66
17B70  
17B81
18G10 
18G35
18G50
53C05
55P62
55R15
55R20 
55S40
55U15
81T70}
\endsubjclass
\keywords
{Twisting cochain, deformation equation, master equation,
models for fiber spaces,
classification of fibrations, formal gauge theory, formal power
series connection, differential graded algebra, differential graded
Lie algebra, homological perturbation theory, cohomological physics}
\endkeywords
\endtopmatter
\document
\rightheadtext{Berikashvili's functor and deformations}

\beginsection Introduction

In \cite\berikash,
N. Berikashvili 
introduced the 
functor $\Cal D$ 
in terms of \lq\lq twisting elements\rq\rq\  or
\lq\lq twisting cochains\rq\rq\  in a differential graded algebra.
At that time 
twisting cochains already had a history in topology and differential
homological algebra,
see e.~g.
\cite\moorefiv\ (where the terminology \lq\lq twisting morphism\rq\rq\ 
is used)
and the historical references to the work of
E. H. Brown and to the S\'eminaire Cartan
in that paper.
With hindsight we see that twisting cochains make precise a certain
piece of structure discovered
by H. Cartan in \cite\cartanon\  and \cite\cartantw\ 
which is behind the notion of transgression.
Berikashvili and his students 
subsequently
studied the functor $\Cal D$ further,
cf. e.~g. \cite{\berikafi, \berikafo, \kadeifou, \sanebfou, \sanebfiv}.
In this note 
we show 
the connections of
this functor 
with  various other developments in the literature, 
in particular with some more
recent ones
related with deformation theory and the master equation.
Our survey
is far from complete at this point, though,
and we apologize for all the relevant citations
which have been omitted.
The formal coincidences between the 
developments we will present here are startling and a
general theory 
unifying and explaining them satisfactorily
as incarnations of 
the same idea
is yet to be found.
Our \lq\lq Leitmotiv\rq\rq\ 
is the observation that some such unification and simplification
is provided by the notion of twisting cochain
and by the ideas underlying deformation theory.
Indeed we will explain below how
Berikashvili's functor $\Cal D$ may be viewed
as a version of the deformation
functor.

I am indebted to J. Stasheff for a number of comments.
It is a pleasure to dedicate this paper to Nodar Berikashvili
for his $70$'th birthday.

\beginsection 1. Berikashvili's functor $D$

The ground ring will be written $R$;
it is a commutative ring with one.
Berikashvili's functor $\Cal D$ may be described briefly in the following way:
Let
$\Cal A$ be a differential graded algebra.
A homogeneous element 
$\tau$ of $\Cal A$
(necessarily of degree $-1$)
is said to be a 
{\it twisting element\/}
or {\it twisting cochain\/}
provided it satisfies the 
equation
$$
D \tau = \tau \tau,
\tag1.1
$$
where $D$ refers to the differential
in $\Cal A$, and where
$\tau \tau$
is the product in $\Cal A$ of
$\tau$ with itself. Denote   the set of twisting elements
by $\Cal T(\Cal A)$.
Next
let $G$ be the group of invertible elements of
$\Cal A^0$, and define an action of $G$
on $\Cal T(\Cal A)$,
to be written
$(x,y) \mapsto x * y,\ x \in G, y \in \Cal T (\Cal A)$,
by means of the formula
$$
x *y = x y x^{-1} + (Dx) x^{-1}.
\tag1.2
$$
This is well defined, that is,
given $x \in G$ and $y \in \Cal T (\Cal A)$,
$x *y$ satisfies the equation (1.1) as well.
This is readily seen by a
straightforward calculation relying on the formulas
$$
(Dx) x^{-1}+x Dx^{-1} = 0, \quad
(Dx^{-1}) x +x^{-1} Dx = 0
$$
which, in turn, follow from $x x^{-1} = 1$.
The set 
of orbits $(\Cal T(\Cal A)) /G$
is written $\Cal D(\Cal A)$
and the assignment to $\Cal A$ of
$\Cal D(\Cal A)$
is a functor from the category of differential graded
algebras to the category of sets.
For later reference, we will say that two twisting elements
$x_1,x_2$ are {\it equivalent\/}
provided they lie in the same $G$-orbit.

Berikashvili applied this functor 
$\Cal D$ in various situations; here is one example
\cite\berikafo: 
Consider a (Serre) fibration
$F \to E \to B$.
Let $\Cal R \roman H_*F$ be a free (additive) resolution 
(over the ground ring)
of the homology $\roman H_*F$ of the fiber
and consider the differential
(trigraded) algebra
$$
\Cal A =C^*(B, \roman{End}(\Cal R \roman H_*F))
$$
with the ordinary cup product as algebra structure.
Given $a \in \Cal A$,
endow $Y=C_* B \otimes (\Cal R \roman H_*F)$
with the twisted differential determined by $a$
and write $Y_a$ for the resulting chain complex.
Berikashvili proved that
the fibration determines an element
$a \in \Cal A$ together with a chain equivalence
$m$ from
$Y_a$
to $C_*E$
satisfying the requisite compatibility properties
in such a way that
$Y_a$ is a model for the
Serre fibration;
he showed that, furthermore,
given an arbitrary element $\overline a$ in the class 
$[a] \in \Cal D(\Cal A)$ of $a \in \Cal T(\Cal A)$,
there is
an isomorphism
from
$Y_{\overline a}$ 
to
$Y_a$
whence there is again
a chain equivalence
$\overline m$ from
$Y_{\overline a}$ to $C_*E$
such that
$Y_{\overline a}$
is a model for the
Serre fibration.
The class of $a$ in
$\Cal T(\Cal A)$
is called the {\it predifferential\/}
of the fibration;
it may be shown that, in a sense,
the predifferential does not depend
on the choice of resolution of the fibre.
The standard filtration
of the model
$Y_a$ 
induces a spectral sequence which is the ordinary
Serre spectral sequence of the fibration.
There is also a cohomology version of the theory.
These models have been used, for example, to study the
section problem for fibrations.
When the homology $\roman H_*F$ of the fiber is free over the ground
ring, the model $Y_a$ comes down to models of the kind introduced
by Hirsch \cite\hirschon.

Prompted by Halperin-Stasheff \cite\halpsttw,
Saneblidze 
\cite\sanebfiv\ 
created a multiplicative version
of the functor $\Cal D$
over the rationals.
To this end,
he reworked
the functor $\Cal D$
with a multiplicative resolution
of the cohomology of the fibre.
He established in particular 
a multiplicative version of
Berikashvili's result, but now
for the cochains of 
a fibre space.
Still in \cite\sanebfiv,
Saneblidze extended this approach so as to incorporate
Hopf algebra structures, and he applied this to 
study the loop space.

Such a multiplicative resolution 
of the cohomology was
already used by
Schlessinger-Stasheff \cite\schlstas\ 
for a similar purpose, that is, to study
deformations of rational homotopy types of spaces
and fibrations;
we will
explain this briefly further below.

{\it Homological perturbation theory\/}
enabled the present author to
construct small free resolutions
by means of appropriate twisting cochains
which,
in turn,
led to complete numerical calculations
in group cohomology \cite\cohomolo,
\cite\modpcoho,
\cite\intecoho\ 
which so far cannot be done
by other methods.

\beginsection 2. The deformation equation

Henceforth we suppose that the ground ring $R$
contains the rationals as a subring.
The equation (1.1)
occurs, perhaps up to a sign, in the literature as
{\it deformation equation\/} or
{\it master equation\/}.
More customary is the corresponding {\it Lie algebra deformation equation\/}
or
{\it Lie algebra master equation\/}:
As an equation on elements of a differential graded Lie algebra
$\Cal L$,
for an element $\Gammat$  of $\Cal L$,
the deformation equation has the form
$$
D \Gammat = \frac 12 [\Gammat, \Gammat] 
\tag2.1
$$
(up to sign).
In particular,
given a cocommutative coaugmented differential graded coalgebra $C$
and a differential graded Lie algebra $\fra g$,
$\Cal L=\roman{Hom}(C,\fra g)$
inherits a differential graded Lie algebra structure
in the standard way,
and 
an element $\tau$ thereof satisfying (2.1), that is to say, 
a morphism $\tau \colon C \to \fra g$
satisfying (2.1), is called a {\it Lie algebra twisting cochain\/},
cf. e.~g. \cite\moorefiv, \cite\quilltwo.
When $\Cal A$ is the universal 
(enveloping)
differential graded algebra
of a 
general
differential graded Lie algebra,
for an element $\Gammat$  of $\Cal L$,
the equations (1.1) and (2.1)
are manifestly equivalent.
Inspection shows that,
given
a differential graded algebra
$\Cal A$ and a solution $\Gammat$ of the equation (1.1),
the operator $d_{\Gammat}$ on $\Cal A$
defined by
$$
d_{\Gammat}
(a) = d a - \Gammat a,\quad
a \in \Cal A,
\tag2.2
$$
yields a new
differential graded algebra structure on $\Cal A$.
Likewise,
given
a differential graded Lie algebra
$\Cal L$ and a solution $\Gammat$ of the equation (2.1),
the operator $d_{\Gammat}$ on $\Cal L$
defined by
$$
d_{\Gammat}
(a) = d a - [\Gammat, a],\quad
a \in \Cal L,
\tag2.3
$$
yields a new
differential graded Lie algebra structure on $\Cal L$.

\beginsection 3. Gauge theory

In view of the well known spectacular results 
obtained by means of gauge theory since the 80's there is
no need to comment on the significance of gauge theory here.
We only reproduce briefly some of the relevant notions
and language
in order to explain how gauge theory relates to Berikashvili's functor.

From now on
the ground ring $R$ is that of the reals or that of the complex numbers.
Consider the space $\Cal A(\xi)$ of connections on a principal bundle
$\xi\colon P \to M$ with structure group $G$ having Lie algebra
$\fra g$. 
Denote by 
$\roman{ad}(\xi)$ the adjoint bundle;
this is the bundle
over $M$ with fiber $\fra g$
which is associated to $\xi$
via the adjoint representation of $G$ on $\fra g$.
The Lie bracket on $\fra g$
induces 
a graded Lie algebra structure
on the $\roman{ad}(\xi)$-valued de Rham forms
$\Omega^*(M,\roman{ad}(\xi))$.
Given a connection $A$ and an $\roman{ad}(\xi)$-valued
1-form $\eta$, the sum $A+\eta$ is again a connection, and its
curvature $F_{A+\eta}$ is given by the well known formula
$$
F_{A+\eta}
= F_A + d_A \eta + \frac 12 [\eta,\eta]
$$
where
$d_A$ denotes the operator of covariant derivative of the connection $A$.
In particular, 
$A+\eta$ has the same curvature as $A$ if and only if
$\eta$ satisfies the equation
$$
d_A \eta + \frac 12 [\eta,\eta] = 0,
\tag3.1
$$
often referred to as the {\it Maurer-Cartan\/}
equation;
notice that this equation for $\eta$
is equivalent
to one of the kind (2.1)
for $\tau=-\eta$.
In particular,
when $A$ is flat,
$(\Omega^*(M,\roman{ad}(\xi)), d_A)$
is a differential graded Lie algebra,
and $A+\eta$ is flat if and only if
$\eta$ satisfies the equation (3.1).
The operator $d_{\tau}$ on
$\Cal L=(\Omega^*(M,\roman{ad}(\xi)))$
given in (2.3),
where $\tau = - \eta$,
then coincides with the operator
$d_{A+\eta}$
of covariant derivative for the flat connection
$A+\eta$.
Moreover, the definition of the operation
(1.2) very much looks like the operation of the group of gauge transformations
on the space of connections.
Thus,
for a general differential graded algebra
$\Cal A$, the value $\Cal D(\Cal A)$
of Berikashvili's functor $\Cal D$
formally looks like a space of gauge equivalence classes of flat
connections.

\medskip\noindent{\bf 4. Chen's formal power series connections}
\smallskip\noindent
A special case of 
twisting element or
solution of the deformation or master equation
was introduced by Chen \cite\chenfou.
Chen's formalism involves formal power series.

We remind the reader that
the ground ring is that of the reals or that of the complex numbers.
In his paper \cite\chenfou,
Chen considers the (bi)graded algebra
$\Omega^*(M,\roman T[[V]])$
of smooth forms on a smooth manifold $M$ with values
in the graded algebra of non-commutative
formal power series $\roman T[[V]]$
on a graded vector space $V$.
In Chen's terminology
\cite\chenfou\ (1.2),
a {\it formal power series connection\/}
is a formal power series in
$\Omega^*(M,\roman T[[V]])$
of (total) degree $-1$
(the degree being appropriately interpreted).
The {\it curvature\/}
of a formal power series connection
$\omega$
is, then, by definition, the element
$$
\kappa = d\omega + \omega \wedge \omega\in \Omega^*(M,\roman T[[V]])
$$
of degree $-2$ where $d$ refers to the ordinary de Rham differential
and $\wedge$ to the cup pairing
induced by the multiplication on
$\roman T[[V]]$;
notice that
$$
\omega \wedge \omega= \frac 12[\omega,\omega] \in
\Omega^*(M,\roman T[[V]])
$$
where $[\ ,\ ]$ denotes the bracket
in
$\Omega^*(M,\roman T[[V]])$
induced from the commutator pairing
$\roman T[[V]] \otimes \roman T[[V]] \to\roman T[[V]]$
in the ordinary way.
In Theorem 1.3.1 of \cite\chenfou,
Chen proves the following:
When $V$ is taken to be the desuspension 
$s^{-1} \widetilde{\roman H}_*(M)$
of the reduced (real) homology
$\widetilde{\roman H}_*(M)$ of $M$,
a splitting
of the ordinary de Rham complex 
$\Omega^*(M)$
into a direct sum of $\roman H^*(M)$,
the coboundaries, and a residual summand
$W$ (written $M$ in Chen's paper),
determines a 
formal power series connection
$\omega \in \Omega^*(M,\roman T[[V]])$
together with a differential
$\partial$ on
$\roman T[[V]]$
which turns the latter into a differential graded algebra in such a way that
$$
\kappa + \omega \partial = 0.
\tag4.1
$$
This formal power series connection 
$\omega$
has a certain uniqueness property.
To reproduce it,
suppose $\omega$ written in the form
$\omega = \omega_1 + \omega_2 + \dots $
where, for $j \geq 1$,
$\omega_j$ is the component whose values lie in
$V^{\otimes j}$.
Now $\omega$
is uniquely determined by the requirement
that
the constituent $\omega_1$ coincide essentially with the projection
of $\Omega^*(M)$ onto $\roman H^*(M)$
(determined by the splitting)
and that,
for $j \geq 2$,
the
$\omega_j$'s be non-zero at most on the residual
summand $W$ of the (chosen) decomposition of $\Omega^*(M)$. 
To explain the significance of
the identity (4.1),
we endow $\Omega^*(M,\roman T[[V]])$
with the differential
$D = d \pm \partial$ 
(appropriate signs being taken)
where $d$ refers to the de Rham
differential and $\partial$ to the operator
on $\Omega^*(M,\roman T[[V]])$
induced by that on $\roman T[[V]]$ denoted by the same symbol.
With this preparation out of the way,
the identity (4.1) says that
the formal power series connection $\omega$
(or rather its negative,
in view of our convention)
is a solution of the deformation or master equation
(1.1), where 
$\Cal A =
(\Omega^*(M,\roman T[[V]]),D)$.
Chen goes on to observe that, in fact,
when $M$ is simply connected,
in view of classical algebraic topology results
related with the cobar construction,
$(\roman T[[V]],\partial)$
is a kind of cobar construction
and hence a model for the real chains of the loop space of $M$.

Chen's construction was
reexamined, generalized, and simplified in the framework
of homological perturbation theory
by Gugenheim in
\cite\gugenhth;
Gugenheim's starting point was the observation that
Chen's splitting
of the ordinary de Rham complex 
quoted above in fact just amounts to a contraction
of chain complexes in the sense of homological perturbation theory;
cf. e.~g. \cite\huebstas.

Suppose that
$V$ is of finite type over the reals
(or complex numbers, according to the case considered).
Appropriately dualized,
Chen's notion of
formal power series connection
then fits into the ordinary twisting cochains framework,
cf. \cite\gulasta.
Indeed,
by means of the tensor coalgebra
$\roman T^{\roman c}[V^*]$ on the dual $V^*$ of $V$,
a formal power series connection may then be written
as a morphism 
$$
\omega
\colon
\roman T^{\roman c}[V^*]
@>>>
\Omega^*(M)
$$
of degree $-1$
with values in the de Rham algebra
$\Omega^*(M)$.
In particular,
when $V$ is taken to be the desuspension 
$s^{-1} \widetilde{\roman H}_*(M)$
of the reduced homology $\widetilde{\roman H}_*(M)$
so that
$V^*$ is the 
suspension 
$s \widetilde{\roman H}^*(M)$
of the reduced cohomology,
the tensor coalgebra
$\roman T^{\roman c}[V^*]$ 
carries an operator
$\delta$ 
compatible with the graded coalgebra structure
(whose dual is the above operator $\partial$ on
$\roman T[[V]]$)
so that
$$
\roman T_{\delta}^{\roman c}[V^*]=(\roman T^{\roman c}[V^*],\delta)
$$
is a model for the bar construction on the de Rham algebra $\Omega^*M$,
that is, a model for the real cochains on the loop space 
of $M$
when $M$ is
simply connnected;
Chen's 
formal power series connection 
satisfying (4.1)
is now precisely a
twisting cochain
$$
\omega
\colon
\roman T_{\delta}^{\roman c}[V^*]
@>>>
\Omega^*(M).
$$
The operator
$\delta$ and formal power series connection
$\omega$ have been constructed in
\cite\gulasta\  by means of
homological perturbation theory.

\medskip\noindent{\bf 5. Deformation theory}\smallskip\noindent
In \cite\haintond, 
Hain and Tondeur interpret Chen's notion of formal power series
connection as the dual of the versal deformation of a trivial connection.
There is indeed an intimate relationship between
Berikashvili's functor, Chen's formal power series connections,
and deformation theory,
in particular 
the deformation theory
of rational homotopy types developed by Schlessinger-Stasheff \cite\schlstas.
We now explain briefly some of the
relevant deformation theory notions.
See e.~g. \cite\gerschtw\ for a recent more detailed account
of deformation theory.

According to a point of view adopted by Deligne and others,
any problem in deformation theory is controlled by a differential graded
Lie algebra, unique up  to homotopy equivalence of differential graded
Lie algebras,
in fact,
up to $L_{\infty}$-equivalence.

Let  $\Cal L$
be a 
differential graded Lie algebra,
the bracket being written
$[\ ,\ ] \colon \Cal L \otimes \Cal L \to \Cal L$.
Recall that
the assignment
$$
\Cal L_0 \times \Cal L_{-1}
@>>>
\Cal L_{-1},
\quad
(\zeta, \alpha) \mapsto
[\zeta, \alpha] - d \zeta
\tag 5.1
$$
yields an action of
$\Cal L_0$
on $\Cal L_{-1}$, the latter being viewed as a vector space,
by infinitesimal
affine transformations.
When $\Cal L_0$ is nilpotent
(or more generally, pro-nilpotent),
the Campbell-Baker-Hausdorff formula turns
$\Cal L_0$
into a (pro-)nilpotent (sometimes Lie) group $\Gamma$
having Lie algebra
$\Cal L_0$,
and the 
$\Cal L_0$-action on $\Cal L_{-1}$
integrates to an 
action
of $\Gamma$ on
$\Cal L_{-1}$
by affine transformations.

Let  $\fra g$
be a 
differential graded Lie algebra,
the bracket being written
$[\ ,\ ] \colon \fra g \otimes \fra g \to \fra g$.
The usual
definition of the functor $\roman{Def}_{\fra g}$
assigns to a local Artinian algebra $\Cal A$
with maximal ideal $\fra m$ the set
$$
\roman{Def}_{\fra g}(\Cal A) =\{\gamma \in \fra g_{-1} \otimes \fra m; 
d \gamma + \frac 12 [\gamma,\gamma] = 0\}\big / \Gamma_{\Cal A};
\tag5.2
$$
here $\Cal L = \fra g \otimes \fra m$
is the induced differential graded $\Cal A$-Lie algebra,
so that the nilpotent Lie algebra $\Cal L_0=\fra g_0 \otimes \fra m$
acts on 
$\Cal L_{-1}=\fra g_{-1} \otimes \fra m$
infinitesimally,
and 
$\Gamma^{\Cal A}$
is the corresponding nilpotent group 
which accordingly acts on
$\Cal L_{-1}$
by affine transformations.
The points of
$\roman{Def}_{\fra g}(\Cal A)$
are the isomorphism classes of the corresponding
{\it Deligne groupoid\/}.
Two 
elements
$\gamma_1,\gamma_2 \in \fra g_{-1} \otimes \fra m$
are said to be {\it equivalent\/}
provided they lie in the same
$\Gamma_{\Cal A}$-orbit.

There is a striking similiarity between
the functor
$\roman{Def}_{\fra g}$
and Berikashvili's functor $\Cal D$.
This similiarity is of course {\it not\/}
a coincidence:
Let $C=\Sigm^{\roman c}]\xi[$
be the
symmetric coalgebra 
with a single {\it cogenerator\/} $\xi$ of degree 0, that is,
the
symmetric coalgebra 
on a free $\Ri$-module with a single basis element
$\xi$ of degree 0.
As an $R$-module, 
$C$  has free generators
$\xi_0 = 1, \xi_1 = \xi,\xi_2, \xi_3, \dots$,
with diagonal map $\Delta$,
counit $\varepsilon$, 
and coaugmentation map $\eta$ being given by
$$
\Delta \xi_n = \sum _{i+j=n} \xi_i \otimes \xi_j,
\quad
\varepsilon (\xi_0) =1, 
\varepsilon (\xi_1) =\varepsilon (\xi_2) = \dots =0, 
\quad
\eta (1) = \xi_0 =1.
$$
Thus, 
when $t\colon C \to R$ denotes the $R$-linear map
which sends $\xi$ to $1$ and is zero on the other generators of $C$,
the ring
$\roman{Hom}(C,R)$  coincides with the complete local ring
$R[[t]]$ of formal power series in the variable $t$,
and the ideal consisting of the formal power series
without constant term is the unique maximal ideal of this ring.
A twisting cochain $\tau \colon C \to \fra g$
is then precisely a formal deformation
in the ordinary sense,
cf. e.~g. \cite\gerschtw\ 
for the notion of formal deformation.
For example when
$\fra g$ arises from the Hochschild complex
of an algebra $B$ in the usual way,
the bracket being the Gerstenhaber bracket
\cite{\gerstenh, \gerschtw},
a twisting cochain $\tau \colon C \to \fra g$
is a formal deformation of $B$.
Furthermore, for a general differential graded
Lie algebra $\fra g$,
$\roman{Hom}(C,\fra g)$
inherits a differential graded Lie algebra structure
as we have already pointed out, and
the coaugmentation map $\eta$ induces 
surjective morphism
$\eta^*$ of differential graded Lie algebras from
$\roman{Hom}(C,\fra g)$
onto
$\fra g$;
thus,
writing
$\Cal L^C$ for the kernel
of $\eta^*$,
we obtain
an extension
$$
0
@>>>
\Cal L^C
@>>>
\roman{Hom}(C,\fra g)
@>{\eta^*}>>
\fra g
@>>>
0
$$
of differential graded Lie algebras.
With respect to the filtration induced from the coaugmentation
filtration of $C$,
the differential graded Lie algebra
$\Cal L^C$ is pronilpotent,
i.~e. an inverse limit of nilpotent Lie algebras.
Thus the infinitesimal action of $\Cal L_0^C$
on 
$\Cal L_{-1}^C$
integrates to an affine action of the corresponding pronilpotent
group $\Gamma^C$
on $\Cal L_{-1}^C$.
The twisting cochains constitute a subset of
$\Cal L_{-1}^C$,
the
$\Gamma^C$-action carries twisting cochains to twisting cochains, and
the value 
$\roman{Def}_{\fra g}(A)$
of the algebra $A=R[[t]]= \roman{Hom}(C,R)$
under the functor
$\roman{Def}_{\fra g}$
is the orbit space of twisting cochains
under the $\Gamma^C$-action;
thus
$\roman{Def}_{\fra g}(A)$
classifies $\fra g$-valued twisting cochains
defined on $C$.
(Technically, one would have to replace the tensor product 
occurring in (5.2) by a completed tensor product when $\fra g_{-1}$ is not
finitely generated.)
On the other hand, we could take
the differential graded algebra
$\Cal A = \roman{Hom}(C, U\fra g)$,
where $U\fra g$ refers to the universal differential
graded algebra of $\fra g$,
and apply
Berikashvili's functor
$\Cal D$ to it;
the resulting object 
$\Cal D(\Cal A)$
again classifies twisting cochains
defined on $C$,
but now with values in all of
$U \fra g$,
and the two notions of equivalence
are somewhat different
since the action of the group
$\Gamma^C$
and that of the group of units $G$ in
$\Cal A^0$
are not in an obvious way related.

Berikashvili's classification procedure,
that is,
that involving the group of units $G$ in
$\Cal A^0$,
provides in fact a classification in terms
of homotopy classes of twisting cochains
with reference to the corresponding notion of homotopy
of twisting cochains: More precisely,
the equation (1.2) says that $x$ provides a homotopy of twisting cochains
between $y$ and $x * y$; see e.~g. \cite\munkholm\ 
for this notion of homotopy of twisting cochains.

The standard formal deformation theory 
classification procedure
involves the group 
corresponding to the degree zero Lie algebra constituent
and the action thereof on
the  Lie algebra constituent
in degree $-1$, though.
Under appropriate circumstances,
an interpretation of the latter classification in terms of
a suitable notion of homotopy
of twisting cochains (more intricate than the naive notion
mentioned above)
may be found in \cite\schlstas.
This raises the question
whether we can isolate the precise relationship between the two
classification procedures.
The whole discussion is of course
valid for {\it any\/}
cocomplete cocommutative differential graded coalgebra 
instead of the present $C$.
Thus we see that {\it Berikashvili's functor may be viewed as a version
of the deformation functor\/}.

We now explain briefly 
some of the links 
between
deformation theory and 
homological perturbation theory.
Consider
a smooth complex manifold $M$.
The differential graded Lie algebra
$\Cal L$
of $\overline \partial$-forms on $M$ 
with values in the holomorphic tangent bundle
$\tau_M$,
the differential being the operator $\overline \partial$,
controls the deformations of the complex structure.
This 
differential graded Lie algebra
is called the {\it Kodaira-Spencer\/} algebra of $M$;
its cohomology 
coincides with the sheaf cohomology 
$\roman H^*(M,\tau_M)$
of $M$ with values
in the holomorphic tangent bundle
$\tau_M$
and plainly inherits a graded Lie algebra structure.
A 1-cochain $\eta$ in $\Cal L$
determines a new complex structure, that is,
a new $\overline \partial$-operator,
if and only $\eta$
satisfies the deformation equation
$$
\overline \partial \eta + \frac 12 [\eta,\eta] = 0,
$$
which, up to sign, is just the  equation (2.1).
Recall that a {\it complex analytic family of complex manifolds over\/} $M$
{\it parametrized by a based open domain\/} $(U,o)$ of $\Bobb C^n$
may be described as a map $\rho \colon U \to \Cal L_{-1}$
whose values lie in 
$\Cal F = \{\eta; \overline\partial \eta + \frac 12 [\eta,\eta] = 0\}$
and which has the requisite properties so that
twisting the complex analytic product $U \times M$ via $\rho$
yields a new 
complex analytic structure on $U \times M$
(the latter here being viewed as the product of two smooth
manifolds only);
in particular, the appropriate adjoint of $\rho$ is a smooth section of the 
smooth vector bundle
on $U \times M$
arising from
pulling back 
the smooth vector bundle
underlying $\Cal L_{-1}$
via the projection $U \times M \to M$, and
$\rho(o) = 0$.
See e.~g. \cite{\kuranone, \kurantwo}
for details.
In coordinates $z_1,\dots,z_n$ on $U$,
$\rho$ may be written as a power series
in these variables with coefficients from
$\Cal L_{-1}$.
This power series being viewed as a formal power series,
the requirement that the values of
$\rho$ lie in
$\Cal F$ is equivalent to 
the negative of $\rho$ being a twisting cochain
$
\Sigm^{\roman c}[V]
@>>>
\Cal L$,
where $V = \roman T_oU \cong\Bobb C^n$
and where
$\Sigm^{\roman c}[V]$
refers to the symmetric coalgebra
on $V$.
Standard techniques reduce the study of such deformations
to cohomology classes
in $\roman H^1(M,\tau_M)$,
and such a class $y$ is unobstructed
if and only if $[y,y]$ is zero
where $[\ ,\ ]$ refers to the induced graded Lie bracket on
$\roman H^*(M,\tau_M)$.
In particular,
according to  Tian and Todorov
\cite{\tianone,\todortwo},
for a Calabi-Yau manifold $M$,
the graded Lie bracket on $\roman H^*(M,\tau_M)$
is zero, and this implies the fact,
first observed by Bogomolov
\cite{\bogomone, \bogomtwo}, 
that the
deformations of a Calabi-Yau manifold are parametrized by an open
subset of
$\roman H^1(M,\tau_M)$.

The $\overline \partial$-forms with values in the 
non-zero exterior powers of
the holomorphic tangent bundle,
endowed with the {\it Fr\"olicher-Nijenhuis bracket\/},
constitute a differential Gerstenhaber algebra
(see e.~g. \cite\banach\ for details).
The Tian-Todorov Lemma says that,
when $M$ is a Calabi-Yau manifold,
this differential Gerstenhaber algebra
has a differential Batalin-Vilkovisky algebra
structure, the requisite generator
being induced from a choice of holomorphic volume form.
In \cite\barakont,
Barannikov and Kontsevich
constructed a formal solution of the {\it master equation\/}
for such a differential Batalin-Vilkovisky algebra.
By means of this formal solution,
they endowed the extended moduli space of complex structures
with a formal Frobenius
manifold structure.
This development, in turn,
was prompted by the mirror conjecture, cf. \cite\wittetwe.
See e.~g. \cite\maninfiv\ 
for more details on Frobenius manifolds.
A purely formal approach
which explains the requisite
differential Batalin-Vilkovisky
algebra structures
in terms of duality for Lie-Rinehart algebras
may be found in \cite\banach\ 
which also contains more references;
some comments about the origin
of Batalin-Vilkovisky
algebras and their relevance in physics will be given in the next section.
In \cite\huebstas,
formal solutions of the master equation
have been obtained much more generally by the methods of homological
perturbation theory.
For related constructions
of twisting cochains
see \cite\gugenmun, 
\cite\perturba\ (2.11),
\cite\munkholm\ (2.2).
The approach in \cite\huebstas\  covers also  the 
requisite twisting cochains
needed for the deformation theory of rational homotopy types
and rational fibrations 
and for the corresponding classification theory
developed by Schlessinger and Stasheff
in \cite\schlstas. In Section 7 below we will briefly
indicate how this relates to the functor $\Cal D$.
For details and history on homological perturbation
theory see e.~g. \cite\huebkade.

\beginsection 6. The master equation in physics

In a series of seminal papers \cite\batviltw, \cite\batvilfo, \cite\batavilk, 
Batalin and Vilkovisky 
studied the quantization of constrained systems
and for that purpose introduced certain differential graded algebras
which have later been named {\it Batalin-Vilkovisky algebras\/}.
Batalin-Vilkovisky algebras have recently become important
in string theory and elsewhere, cf. e.~g.
\cite\barakont,
\cite\bv,
\cite\huebstas,
\cite\maninfiv,
\cite\stashnin.
String theory leads to the
mysterious mirror conjecture,
see e.~g. \cite\wittetwe.
In this context, 
the
quantum Batalin-Vilkovisky master equation
has the form of the Maurer-Cartan equation
for a flat connection
while the classical version has the form of the
integrability condition of deformation theory.
See \cite\stashnin\  for more details and references.

\beginsection 7. Concluding remarks

We have seen 
that Berikashvili's functor may be viewed as a version
of the deformation functor, and that (formal) deformations
may be viewed as twisting cochains defined on
coaugmented cocommutative differential graded coalgebras.
Chen's formal power series connections
are also subsumed under the notion of twisting cochain,
and the curvature of a Chen formal power series connection
being zero amounts to the 
corresponding master equation being satisfied.
Thus, when twisting cochains defined on not necessarily
cocommutative coalgebras are 
viewed as non-commutative generalizations of deformations,
we can perhaps give a meaning to the 
point of view of Hain and Tondeur 
that Chen's construction of formal power series
connection should be
viewed as  the versal deformation of something
(or the dual thereof).

In the same rite,
the deformation theory for rational homotopy types
developed by Schlessinger and Stasheff \cite\schlstas\ 
involves as an ingredient the (filtered multiplicative) Halperin-Stasheff
model \cite\halpsttw\ for a differential graded commutative algebra
and in particular the master equation.
Saneblidze's multiplicative version of
the multiplicative
functor $\Cal D$ \cite\sanebfiv\ 
also relies on 
the Halperin-Stasheff
model and the requisite master equation.
According to
Saneblidze's 
Theorem 4.8
in \cite\sanebfiv, 
under the circumstances of that theorem,
for appropriate spaces $X$, 
the object
$\Cal D(X)$ introduced there (defined in terms of the
corresponding Sullivan-de Rham complex)
classifies fiber homotopy types over $X$ 
where
the homotopy type of the fiber is fixed in advance,
while Schlessinger-Stasheff 
obtain  such a classification in Theorem 9.6 of \cite\schlstas\,
in terms of their moduli space arising from
their differential graded Lie algebra.
Thus the two classification theories are plainly related,
but can we make the 
relationship precise and explicit?
Another intriguing observation,
which may be found just before Theorem 5.5 in \cite\sanebfiv,
is that
there the requisite
\lq\lq Hopf predifferential\rq\rq\ 
may be obtained as a Chen formal power series
connection,
in view of what is said in \cite\hainone.
Thus, what is the precise relationship
between the multiplicative functor $\Cal D$ and
Chen's formal power series connections?

Perhaps some of these questions
related
with Berikashvili's functor
will be answered in the future.

\bigskip 
\widestnumber\key{999}
\centerline{References}
\smallskip\noindent

\ref \no \barakont
\by S. Barannikov and M. Kontsevich
\paper Frobenius manifolds and formality of Lie algebras of polyvector fields
\paperinfo {\tt alg-geom/9710032}
\jour Internat. Res. Notices
\vol 4
\yr 1998
\pages 201--215
\endref

\ref \no \batviltw
\by I. A. Batalin and G. S. Vilkovisky
\paper Quantization of gauge theories
with linearly dependent generators
\jour  Phys. Rev. 
\vol D 28
\yr 1983
\pages  2567--2582
\endref

\ref \no \batvilfo
\by I. A. Batalin and G. S. Vilkovisky
\paper Closure of the gauge algebra, generalized Lie equations
and Feynman rules
\jour  Nucl. Phys. B
\vol 234
\yr 1984
\pages  106-124
\endref

\ref \no \batavilk
\by I. A. Batalin and G. S. Vilkovisky
\paper Existence theorem for gauge algebra
\jour Jour. Math. Phys.
\vol 26
\yr 1985
\pages  172--184
\endref

\ref \no \berikash
\by N.A. Berikashvili
\paper On the differentials of a spectral sequence
\jour Bulletin of the Academy of Sciences of the Georgian SSR
\vol 51
\yr 1968
\pages 9--14
\finalinfo  (Russian)
\endref

\ref \no \berikafi
\by N. A. Berikashvili
\paper On the differentials of spectral sequences
\jour Trudy Tbilisi Math. Institute
\vol 51
\yr 1976
\pages 1-105
\finalinfo  (Russian) 
\endref

\ref \no \berikafo
\by N. A. Berikashvili
\paper On the obstruction theory in fiber spaces
\jour Bulletin of the Academy of Sciences of the Georgian SSR
\vol 125
\yr 1987
\pages 473--475
\finalinfo  (Russian) 
\endref

\ref \no \bogomone
\by F. A. Bogomolov
\paper Hamiltonian K\"ahler varieties
\jour Sov. Math. Dokl.
\vol 19
\yr 1978
\pages 1462--1465
\finalinfo
Translated from: Dokl. Akad. Nauk SSR 243 (1978) 1101--1104
\endref

\ref \no \bogomtwo
\by F. A. Bogomolov
\paper K\"ahler manifolds with trivial canonical class
\paperinfo Preprint, Institut des Hautes Etudes
Scientifiques
\yr 1981
\pages 1--32
\endref

\ref \no \cartanon
\by H. Cartan
\paper Notions d'alg\`ebre diff\'erentielle;
applications aux groupes de Lie et aux vari\'et\'es
o\`u op\`ere un groupe
de Lie
\jour Coll. Topologie Alg\'ebrique
\paperinfo Bruxelles
\yr 1950
\pages  15--28
\endref

\ref \no \cartantw
\by H. Cartan
\paper La transgression dans un groupe de Lie et dans un espace
fibr\'e principal
\jour Coll. Topologie Alg\'ebrique
\paperinfo Brussels
\yr 1950
\pages  57--72
\endref

\ref \no \cartanei
\by H. Cartan and S. Eilenberg
\book Homological Algebra
\publ Princeton University Press
\publaddr Princeton
\yr 1956
\endref

\ref \no \chenfou
\by K. T. Chen
\paper Extension of $C^{\infty}$
Function Algebra by Integrals and Malcev Completion of $\pi_1$
\jour Advances in Mathematics
\vol 23
\yr 1977
\pages 181--210
\endref

\ref \no \eilmactw
\by S. Eilenberg and S. Mac Lane
\paper On the groups ${\roman H(\pi,n)}$. I.
\jour Ann. of Math.
\vol 58
\yr 1953
\pages  55--106
\moreref
\paper II. Methods of computation
\jour Ann. of Math.
\vol 60
\yr 1954
\pages  49--139
\endref

\ref \no \gerstenh
\by M. Gerstenhaber
\paper On the deformation of rings and algebras.I.
\jour Ann. of Math.
\vol 79
\yr 1964
\pages  59--103
\moreref
\paper II.
\jour Ann. of Math.
\vol 84
\yr 1966
\pages  1--19
\moreref
\paper III.
\jour Ann. of Math.
\vol 88
\yr 1968
\pages  1--34
\moreref
\paper IV.
\jour Ann. of Math.
\vol 99
\yr 1974
\pages  257--276
\endref

\ref \no \gerschtw           
\by M. Gerstenhaber and Samuel D. Schack
\paper Algebraic cohomology and deformation theory
\paperinfo In: Deformation theory of algebras and structures
and applications, M. Hazewinkel and M. Gerstenhaber, eds.
\pages 11--264
\publ Kluwer
\yr 1988
\publaddr Dordrecht
\endref

\ref \no \gugenhtw
\by V.K.A.M. Gugenheim
\paper On the chain complex of a fibration
\jour Illinois J. of Mathematics
\vol 16
\yr 1972
\pages 398--414
\endref

\ref \no \gugenhth
\by V.K.A.M. Gugenheim
\paper On a perturbation theory for the homology of the loop space
\jour J. of Pure and Applied Algebra
\vol 25
\yr 1982
\pages 197--205
\endref

\ref \no \gugenlam
\by V.K.A.M. Gugenheim and L. Lambe
\paper Perturbation in differential homological algebra
\jour Illinois J. of Mathematics
\vol 33
\yr 1989
\pages 566--582
\endref

\ref \no \gulasta
\by V.K.A.M. Gugenheim, L. Lambe, and J.D. Stasheff
\paper Algebraic aspects of Chen's twisting cochains
\jour Illinois J. of Math.
\vol 34
\yr 1990
\pages 485--502
\endref

\ref \no \gulstatw
\by V.K.A.M. Gugenheim, L. Lambe, and J.D. Stasheff
\paper Perturbation theory in differential homological algebra. II.
\jour Illinois J. of Math.
\vol 35
\yr 1991
\pages 357--373
\endref

\ref \no \gugenmun
\by V.K.A.M. Gugenheim and H. J. Munkholm
\paper On the extended functoriality of Tor and Cotor
\jour J. of Pure and Applied Algebra
\vol 4
\yr 1974
\pages  9--29
\endref
\ref \no \hainone
\by R. H. Hain
\paper Twisting cochains and duality between minimal algebras
and minimal Lie algebras
\jour Trans. Amer. Math. Soc.
\vol 277
\pages 397--411
\yr 1983
\endref

\ref \no \haintond
\by R. H. Hain and P. Tondeur
\paper The life and work of Kuo-Tsai Chen
\jour Illinois J. of Math.
\vol 34
\yr 1990
\pages 
\endref

\ref \no \halpstas
\by S. Halperin and J.D. Stasheff
\paper Differential algebra in its own rite
\paperinfo Proc. Adv. Study Inst. Alg. Top.,
August 10--23, 1970, Aarhus, Denmark
\pages 567--577
\endref

\ref \no \halpsttw
\by S. Halperin and J.D. Stasheff
\paper Obstructions to homotopy equivalences
\jour Advances in Math.
\vol 32
\yr 1979
\pages 233--278
\endref

\ref \no \hirschon
\by G. Hirsch
\paper Sur les groupes d'homologie des espaces fibr\'es
\jour Bull. Soc. Math. Belgique
\vol 6
\yr 1953
\pages 79--96
\endref

\ref \no \perturba
\by J. Huebschmann
\paper Perturbation theory and free resolutions for nilpotent
groups of class 2
\jour J. of Algebra
\yr 1989
\vol 126
\pages 348--399
\endref

\ref \no \cohomolo
\by J. Huebschmann
\paper Cohomology of nilpotent groups of class 2
\jour J. of Algebra
\yr 1989
\vol 126
\pages 400--450
\endref

\ref \no \modpcoho
\by J. Huebschmann
\paper The mod $p$ cohomology rings of metacyclic groups
\jour J. of Pure and Applied Algebra
\vol 60
\yr 1989
\pages 53--105
\endref

\ref \no \intecoho
\by J. Huebschmann
\paper Cohomology of metacyclic groups
\jour Trans. Amer. Math. Soc.
\vol 328
\yr 1991
\pages 1-72
\endref
\ref \no \bv
\by J. Huebschmann
\paper Lie-Rinehart algebras, Gerstenhaber algebras, and Batalin-
Vilkovisky algebras
\jour Annales de l'Institut Fourier
\vol 48
\yr 1998
\pages 425--440
\endref

\ref \no \banach
\by J. Huebschmann
\paper Differential Batalin-Vilkovisky algebras arising from
twilled Lie-Rinehart algebras
\paperinfo Proceedings of the Poissonfest, Warsaw, 1998
\jour Banach Center publications (to appear)
\yr 1999
\endref

\ref \no \huebkade
\by J. Huebschmann and T. Kadeishvili
\paper Small models for chain algebras
\jour Math. Z.
\vol 207
\yr 1991
\pages 245--280
\endref

\ref \no \huebstas
\by J. Huebschmann and J. D. Stasheff
\paper Formal solution of the master equation via HPT and
deformation theory
\paperinfo Preprint June 1999
\endref

\ref \no \kadeifou
\by T. Kadeishvili
\paper The predifferential of a twisted product
\jour Russian Math. Surveys
\vol 41
\yr 1986
\pages 135--147
\endref

\ref \no \kosmathr
\by Y. Kosmann-Schwarzbach 
\paper Exact Gerstenhaber algebras and Lie bialgebroids
\jour  Acta Applicandae Mathematicae
\vol 41
\yr 1995
\pages 153--165
\endref

\ref \no \kuranone
\by M. Kuranishi
\paper On the locally complete families of ccmplex analytic structures
\jour Ann. of Math.
\vol 75
\yr 1962
\pages 536--577
\endref

\ref \no \kurantwo
\by M. Kuranishi
\book Deformations of compact complex manifolds
\publ Les presses de l'Universit\'e de Montr\'eal
\publaddr Montr\'eal
\yr 1971
\endref

\ref \no \maninfiv
\by Yu. I. Manin
\paper Three constructions of Frobenius manifolds
\paperinfo {\tt math.QA/9801006}
\jour Atiyah-Festschrift (to appear)
\endref

\ref \no \moorefiv
\by J. C. Moore
\paper Differential homological algebra
\jour Actes, Congres intern. math. Nice
\publ Gauthiers-Villars
\publaddr Paris, 1971
\yr 1970
\pages 335--339
\endref

\ref \no \munkholm
\by H. J. Munkholm
\paper The Eilenberg--Moore spectral sequence and strongly homotopy
multiplicative maps
\jour J. of Pure and Applied Algebra
\vol 9
\yr 1976
\pages  1--50
\endref

\ref \no \quilltwo
\by D. Quillen
\paper Rational homotopy theory
\jour Ann. of Math.
\vol 90
\yr 1969
\pages  205--295
\endref

\ref \no \sanebfou
\by S. Saneblidze
\paper Obstructions to the section problem in fibre bundles
\jour manuscripta math.
\vol 81
\yr 1993
\pages 95--111
\endref

\ref \no \sanebfiv
\by S. Saneblidze
\paper The set of multiplicative predifferentials and the rational 
cohomology algebra of a fiber space
\jour J. of Pure and Applied Algebra
\vol 73
\yr 1991
\pages 277--306
\endref

\ref \no \schlstas
\by M. Schlessinger and J. Stasheff
\paper Deformation theory and rational homotopy type
\jour Pub. Math. Sci. IHES (to appear)
\paperinfo preprint 1979;  new version July 13, 1998
\endref

\ref \no \stashnin
\by J. D. Stasheff
\paper Deformation theory and the Batalin-Vilkovisky 
master equation
\paperinfo in: Deformation Theory and Symplectic Geometry,
Proceedings of the Ascona meeting, June 1996,
D. Sternheimer, J. Rawnsley, S.  Gutt, eds.,
Mathematical Physics Studies, Vol. 20
\publ Kluwer Academic Publishers
\publaddr Dordrecht/Boston/London
\yr 1997
\pages 271--284
\endref

\ref \no \tianone
\by G. Tian
\paper A note on K\"ahler manifolds with $c_1=0$
\paperinfo preprint
\endref

\ref \no \todortwo
\by A. N. Todorov
\paper The Weil-Petersson geometry
of the moduli space of $\fra {su}(n \geq 3)$ (Calabi-Yau)
manifolds, I.
\jour Comm. Math. Phys.
\vol 126
\yr 1989
\pages 325--346
\endref

\ref \no \wittetwe
\by E. Witten
\paper Mirror manifolds and topological field theory
\paperinfo in: Essays on mirror manifolds, 
S. T. Yau, ed.
\publ International Press Co.
\publaddr Hong Kong
\yr 1992
\pages 230--310 
\endref
\enddocument